\documentclass[a4paper,10pt,french]{article}
\usepackage{amsfonts}
\usepackage[T1]{fontenc}
\usepackage[ansinew]{inputenc}
\usepackage{babel}

\usepackage{graphicx}

\usepackage[top=2.5cm,bottom=2.5cm,right=2cm,left=2cm]{geometry}

\begin{document}

\font\cm=cmr10 at 7.5pt

\font\cm=cmr10 at 7.5pt

\font\cmdpp=cmdunh10 scaled 1100

\font\tf=cmssbx10 scaled 1400

\font\tfp=cmssbx10 scaled 1200

\font\tfpp=cmssbx10 scaled 1000

\font\tfppp=cmssbx10 scaled 950

\font\zito=cmssbx10 scaled 1000

\font\titrefontp=cmssbx10 scaled 1400

\font\titrefontppp=cmssbx10 scaled 1000

\font\cmd=cmdunh10 scaled 1500

\font\cmddd=cmdunh10 scaled 2000

\font\cmdp=cmdunh10 scaled 1300

\font\cmdpp=cmdunh10 scaled 1100

\font\tata=cmcsc10 at 10pt

\font\titrefontpp=cmssbx10 scaled 1200

\baselineskip 13.7pt

\font\zazou = cmb10 at 16pt

\font\ib = cmbxsl10

\font\sc=cmcsc10 at 10pt

\font\scc=cmcsc10 at 13pt

\font\pg=cmss10 scaled 1000

\font\pgp=cmss10 scaled 900

\font\itp=cmsl10 scaled 900

\font\pgg=cmss10 scaled 1300

\font\pggp=cmss10 scaled 1200

\font\tfg=cmssbx10 scaled 1800

\font\tf=cmssbx10 scaled 1400

\font\tfp=cmssbx10 scaled 1200

\font\tfpp=cmssbx10 scaled 1000

\font\ib = cmbxsl09

\font\matha=msam10 scaled 1000

\font\mathb=msbm10 scaled 1000
\def\mathbbb#1{\hbox{\mathb #1\relax}}
\let\Bbb=\mathbbb

\centerline{\zazou Area and perimeter foliations on spaces of polygons}

\smallskip
\centerline{by}
\smallskip

\centerline{{Aziz EL KACIMI ALAOUI \& Abdellatif ZEGGAR}}
\smallskip

\centerline{(February 2019)}

\vskip0.5cm

\medskip

\noindent {\parindent=1.2cm\narrower {\pg We describe all families of star-shaped $n$-polygons
in the Euclidean plane with prescribed perimeter and area;
they are leaves of a foliation ${\cal F}^{\;(\ast )}$ on the space ${\cal P}_n^\ast $ of
star-shaped polygons. By the way, we study some
geometric properties of convex polygons, for instance their inscriptibility in a circle
and their regularity in relation with the perimeter and the area. \par }}

\vskip0.8cm

\noindent \noindent We denote by ${\Bbb E}$  an Euclidean vector plane equipped with
its canonical affine structure and an orientation given by an orthonormal basis. The word `isometry' means a
transformation of ${\Bbb E}$ which preserves
the distance; necessarily it is an affine transformation of ${\Bbb E}$.
The origin of ${\Bbb E}$ is denoted $O$. For $A, B\in {\Bbb E}$ we denote $\overrightarrow{AB}$ the vector $B-A$.
Let $n\geq 3$ be a natural integer.

\medskip

{\tfp 0. Introduction}

\medskip

\noindent A {\it figure} of ${\Bbb E}$ is just a part of it. But usually this name is given only to a one having
a certain peculiarity: we see it all (it is bounded) or at least we understand how it is  made to guess
its behavior when it escapes our view, like a straight line... and its outline has a little regularity.
A {\it polygon} is an example of such figure: it is bounded and bordered by a finite number of segments called  {\it sides}
or {\it edges}. To each polygon, one can associate invariants, among which are two real numbers
that play an important role: the {\it perimeter} and the {\it area}.

\smallskip

The polygons of the plane are numerous and their shapes and sizes are varied. So the question of
their equivalence therefore arises naturally. But in which sense?

\smallskip

From the set point of view, two polygons $\wp $ and $\wp'$ are always equivalent: there is a bijection of ${\Bbb E}$
which sends one on the other. But as we are in a Euclidean plane,
we would like this bijection to preserve
properties related to the affine Euclidean structure. There are
several notions of equivalence; here are some of them (those that
will interest us directly). We will say that $\wp $ and $\wp'$ are:

\begin{enumerate}

\item {\it isometric} if there exists an isometry $f: {\Bbb E}\longrightarrow {\Bbb E}$ such that $f(\wp )= \wp'$.
We can superimpose them; and we can still do that without going out
of the plan if $\wp $ and $\wp'$ are {\it directly isometric}, that is $f$ preserves the orientation;

\item {\it similar}  if there is a similarity $f: {\Bbb E}\longrightarrow {\Bbb E}$ such that $f(\wp )= \wp'$.
In a way, one of them is an enlargement of the other (as for the photos);

\item {\it equivalent} (shortly) if they have the same area. In this case, we can always go from one to another by cutting and geometric gluing;

\item {\it isoperimetric} if they have the same perimeter.

\end{enumerate}

The isometric equivalence is the strongest and implies  all the others. So it is too rigid to be `useful':
two isometric polygons differ only in the positions they occupy in the plane.
It is rather the equivalences 3 and 4 that will occupy us here.

\smallskip

This paper takes its origin from the following question ($\ast $):  {\it Are there two non-isometric triangles with same perimeter and same area}?
It led us first to the construction of a foliation ${\cal F}^{\;(\ast )}$  on the space of triangles: each leaf consists of triangles
having same perimeter and same area. Then we made a more general study for the space of star-shaped polygons; this space contains convex
polygons for which some properties related to area and perimeter have also been studied. Some of them are certainly known, but little absent in
the `geometric literature', which motivated us to insert them in this text.

\vskip0.4cm

\hrule

\vskip0.4cm

\noindent ($\ast $) It was asked by Geoffrey Letellier
to his teacher Valerio Vassallo who put it in turn to the first author. The latter has built for this purpose the foliation ${\cal F}$
on the space of the triangles. But just to the question asked, Geoffrey
himself responded by constructing a one parameter family of isosceles triangles having the same perimeter and the same area (his
example is given in subsection 4.6).

\medskip

{\ib The integer $n$ we consider in all this paper will be greater than or equal to $3$.}

\medskip
{\tfp 1. Preliminaries}

\medskip

\noindent {\tfpp 1.1. Definition.} {\it A non degenerate  {\it polygon} of ${\Bbb E}$ is an element $\wp =(M_1,\cdots ,M_n)$ of ${\Bbb E}^n$ such that}:

(i) {\it for $i\neq j$ the point $M_i$ is distinct from $M_j$};

(ii) {\it for any $k\in \{1,\cdots ,n\}$, the oriented angle $\widehat{M_k}=(\overrightarrow{M_kM_{k+1}},\overrightarrow{M_kM_{k-1}})$ has its measure in
$]0,2\pi [\setminus \{ \pi\} $.}

\smallskip

By convention, we set $M_0=M_n$ and $M_{n+1}=M_1$. The orientation of the angle $\widehat{M_k}$ is the same as the orientation of the  trigonometric
circle  centered at the point $M_k$.

\smallskip

The points $M_k$ and the segments $[M_kM_{k+1}]$  are respectively the
{\it vertices} and the {\it sides} of the polygon $\wp $. If $M_i$ and $M_j$  are two non successive vertices, that is  $\vert i-j\vert > 1$,
we say that the segment  $M_iM_j$ is a {\it diagonal} of the  polygon.

\smallskip

Recall that a polygon is said to be:

\smallskip

- {\it equilateral} if all its  sides have the same length;

- {\it inscribable} if all its vertices are on a same circle;

- {\it regular} if it is both inscribable and equilateral.

\smallskip

The set of all $n$-polygons of the plane ${\Bbb E}$ will be denoted
$\widetilde{\cal{P}}_n$.
Let  $f$ be an isometry of  the Euclidean plane ${\Bbb E}$.  The  image   $\wp'=(M'_1,\cdots ,M'_n)=\left(f(M_1),\cdots ,f(M_n)\right)$ of a polygon
$\wp =(M_1,\cdots ,M_n)$
is a polygon of  ${\Bbb E}$ such that:
$$\widehat{M'_k}=\widehat{M_k}\hskip0.2cm\hbox{and}\hskip0.2cm\vert \vert \overrightarrow{M'_kM'_\ell }
\vert \vert =\vert \vert \overrightarrow{M_kM_\ell }\vert \vert  \; \; \hbox{for}\; \; (k,\ell )\in \{1,\cdots ,n\}^2.
\leqno{(1.1)}$$
So we have a natural action:
$$\hbox{Isom}({\Bbb E})\times \widetilde{\cal{P}}_n\longrightarrow \widetilde{\cal{P}}_n, \; \big( f,(M_1,\cdots ,M_n)\big) \longmapsto
\big( f(M_1),\cdots ,f(M_n)\big) \leqno{(1.2)}$$
where $\hbox{Isom}({\Bbb E})$ is  the group of the affine isometries  of  ${\Bbb E}$. The quotient space of this action will
be denoted:
$${\cal P}_n=\widetilde{\cal P}_n/\hbox{Isom}({\Bbb E}).\leqno{(1.3)}$$

 The elements of ${\cal P}_n$  are called
{\it geometric polygons}  of ${\Bbb E}$.
\smallskip

A geometric polygon of ${\Bbb E}$ is said to be {\it equilateral} (resp. {\it inscribable}) if it admits a representative
which is equilateral (resp. inscribable).
\smallskip

For an element $\wp =(M_1,\cdots ,M_n)$ of $\widetilde{\cal P}_n$, we will use the following notations:
$$\cases{<M_1,\cdots ,M_n>  \hbox{ is the equivalence class of
$(M_1,\cdots ,M_n)$  in $\widetilde{\cal P}_n$,}\cr {}\cr
x_k=M_kM_{k+1}=\vert \vert \overrightarrow{M_kM_{k+1}}\vert \vert  \; \;
\hbox{for}\; \; k\in \{1,\cdots ,n-2\},\cr {}\cr
t_k=M_nM_k=\vert \vert \overrightarrow{M_nM_k}\vert \vert  \; \; \hbox{for}\; \;
k\in \{1,\cdots ,n-1\} .
}\leqno{(1.4)}$$
The positive numbers $t_1,x_1,\cdots ,x_{n-2},t_{n-1}$
are the lengths  of the sides  and $t_2,t_3,\cdots ,
t_{n-2}$ are the lengths of the diagonals from the vertex $
M_n$ (see the picture bellow for the case of the hexagon). We have $(n-2)$ lengths of  type $x_k
$ and $(n-1)$ lengths of type $t_k$.
\smallskip

Since an isometry preserves the distance in the Euclidean plane ${\Bbb E}$,
$(t_1,x_1,t_2,x_2,\cdots ,t_{n-2},x_{n-2},t_{n-1})$ does not depend on the choice of the representative
 $(M_1,\cdots ,M_n)$ of the geometric polygon
$<M_1,\cdots ,M_n>$.

$$\includegraphics[height=4.5cm]{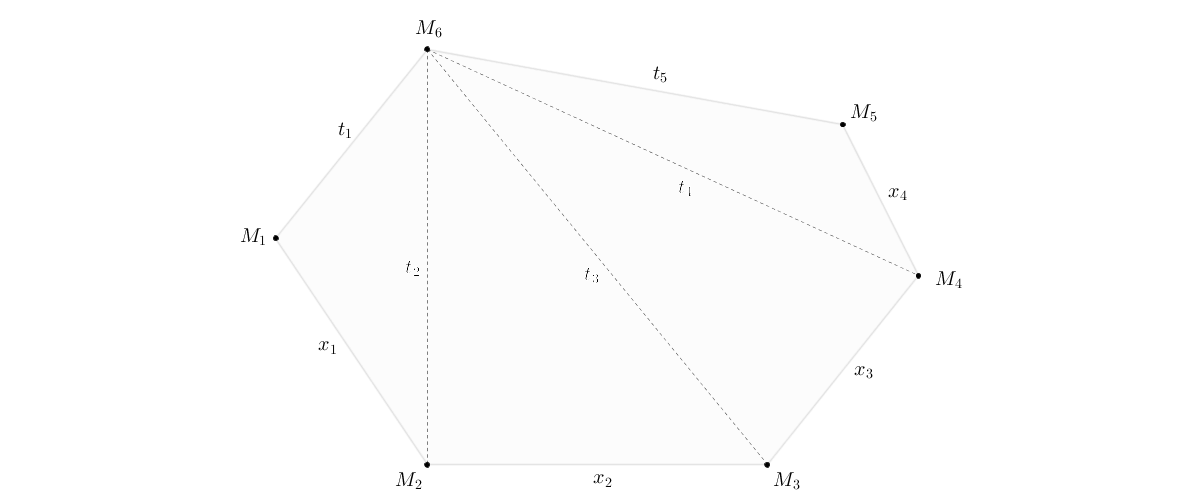}$$

For any  triangle $(M_k,M_n,M_{k+1})$ with $
1\leq k\leq n-2$, the lengths $t_k,x_k, t_{k+1}$ satisfy the inequalities:
$$\cases{0<t_k<x_k+t_{k+1}\cr
0<x_k<t_k+t_{k+1}\cr
0<t_{k+1}<x_k+t_k}\leqno{(1.5)}$$
that is, $(t_k,x_k,t_{k+1})$
is an element of the open set ${\cal V}$ of ${\Bbb R}^3$
consisting of the triplets $(x,y,z)$ satisfying:
$$\cases{
0<x<y+z\cr
0<y<x+z \cr
0<z<x+y.} $$
Moreover, one can verify by induction that, for the lengths
of the sides of a polygon, each length is strictly smaller than the sum of the others.
For instance, we have the following nice:

\smallskip

\noindent {\tfpp 1.2. Theorem} [Pen]. {\it For any natural integer $n\geq
3$ and any $n$-tuple $u=(u_1,\cdots ,u_n)$ of positive real numbers, there exists a unique inscribable geometric  polygon $<M_1,\cdots ,M_n>$ such that
$M_kM_{k+1}=u_k$ for $k\in \{1,\cdots , n\} $ if and only if, for any $j\in \{ 1,\cdots ,n\} $, we have the inequality: $u_j<\displaystyle \sum _{k\neq j}u_k.$}

\smallskip

\noindent {\tfpp 1.3. Definition.} {\it A  $n$-polygon  $\wp =(M_1,\cdots ,M_n)$ of ${\Bbb E}$  is {\ib convex} if, for any
$k\in \{1,\cdots ,n\}$, the oriented angle $\widehat{M_k}=(\overrightarrow{M_kM_{k+1}},\overrightarrow{M_kM_{k-1}})$ has its measure in
$]0,\pi [$. It is called {\ib  star-shaped polygon} with respect to the vertex $M$ if, for any vertex $N\in \{ M_1,\cdots ,M_n\} $ the open segment
$]MN[$ is contained in the interior of $\wp $. Of course, any convex polygon is a  star-shaped polygon with respect to any of its vertices.}

\smallskip

Star-shaped polygons $(M_1,\cdots ,M_n)$ with respect to the vertex $M_n$ form a subspace $\widetilde{\cal P}_n^\ast $ of $\widetilde{\cal P}_n$, invariant
under the action of $\hbox{Isom}({\Bbb E})$ (on $\widetilde{\cal P}_n$).

One can easily see that any element of $\widetilde{\cal P}_n^\ast $ is isometric to a unique  star-shaped polygon
$\wp = (M_1,\cdots ,M_n)$ with respect to $M_n$ such that:
$$\cases{M_n=M_0=O \hbox{ the origin of ${\Bbb E}$,}\cr {}\cr  M_1 \hbox{ has coordinates $(t_1,0)$ with $t_1>0$,}\cr  {}\cr \hbox{for any $k\in \{1,\cdots ,n-2\}$,
the measure of the angle $\widehat{M_kOM_{k+1}}=(\overrightarrow{OM_k},\overrightarrow{OM_{k+1}})$ is in $]0,\pi [$.}}\leqno{(1.6)}$$

\smallskip

${\cal P}_n^\ast $ will denote the space of the star-shaped polygons satisfying these conditions.

\smallskip

$$\includegraphics[height=6.5cm]{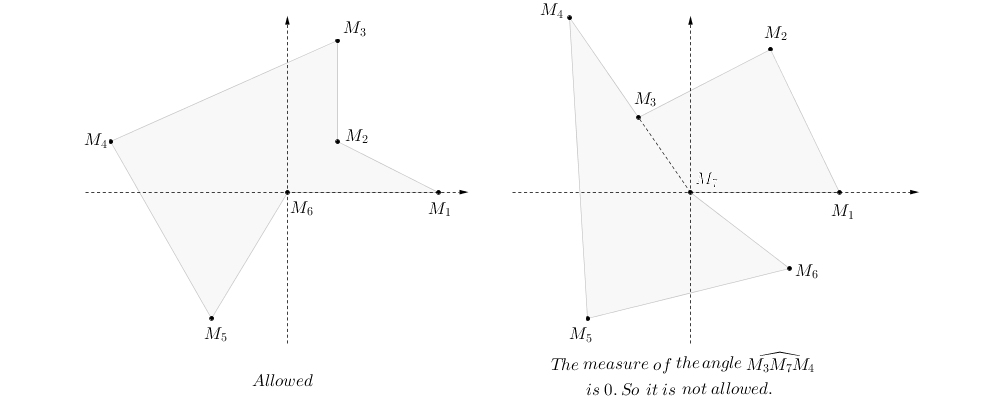}$$

\smallskip

In the following subsection we will recall some metric properties in a triangle which will be very useful for determining the space  of polygons on which we will
define area and perimeter foliations in section 3.

\smallskip

\noindent {\tfpp 1.4.} Let $OMN$ be a non degenerate triangle ; we set $OM=t$, $ON=s$, $MN=x$ and
$\alpha =\widehat{MON}=(\overrightarrow{OM},\overrightarrow{ON})$. All the
real numbers $t$, $s$, $x$ are positive and $\alpha \in ]0,\pi [$.  For the triangle $OMN$ we denote $r>0$ the radius of
its inscribed circle, $p$ its perimeter and $a$ its area. We have the following
well known formulaes:
$$p=t+x+s, \hskip1cm a={1\over 4}\sqrt{(t+x+s)(t-x+s)(t+x -s)(-t+x+s)}\hskip0.3cm \hbox{(H\'eron's formula)}$$
and:
$$r={a\over p}, \hskip1cm \sin (\alpha )={x\over {2r}}={{px}\over {2a}}=
{{2x(t+x+s)}\over {\sqrt{(t+x+s)(t-x+s)(t+x -s)(-t+x+s)}}}.$$

For a given $d\in ]0,1[$ there are two numbers $\alpha'(d)\in ]0,{\pi \over 2}[$ and
$\alpha''(d)\in ]{\pi \over 2},\pi [$ such that $\sin (\alpha'(d))=\sin (\alpha''(d))=d$.
We take $d={{px}\over {2a}}$ (which depends on $(t,x,s)$) and
we define the function $\alpha :(]0,+\infty [)^3\longrightarrow ]0,\pi [$ by:
$$\alpha (t,x,s)=\cases{\alpha'(d(t,x,s))&\hbox{ if $x^2<t^2+s^2$}\cr {}\cr
{\pi \over 2} & \hbox{ if $x^2=t^2+s^2$}\cr {}\cr
\alpha''(d(t,x,s))&\hbox{ if $x^2>t^2+s^2$.}}$$
It is well defined and continuous.

\smallskip

Now consider
the open set $\Omega _n$ of $({\Bbb R}_+^\ast )^{2n-3}$ given by:
$$\Omega _n=\left\{(t_1,x_1,\cdots ,t_{n-2},x_{n-2},t_{n-1}): k=1,\cdots ,n-2,\; (t_k,x_k,t_{k+1})\in {\cal V}, \; \; \sum_{k=1}^{n-2}\alpha (t_k,x_k,t_{k+1})<2\pi \right\} .\leqno{(1.7)}$$

\smallskip

Let ${\cal L}_n:{\cal P}_n^\ast \longrightarrow \Omega _n$ be the map
 ${\cal L}_n(\wp )={\cal L}_n(<M_1,\cdots ,M_n>)=(t_1,x_1,t_2,x_2,\cdots ,t_{n-2},x_{n-2},t_{n-1})$
where the real numbers $t_1,x_1,t_2,x_2,\cdots ,t_{n-2},x_{n-2},t_{n-1}$ are given by the formulaes (1.4).

\smallskip

From now on we will identify the geometric  star-shaped polygons to the points of the open set $\Omega _n$ by the map ${\cal L}_n$.
This identification ${\cal P}_n^\ast \simeq \Omega _n$
will enables one to study easily some properties of the space of geometric star-shaped polygons.

\smallskip

Let $\varphi : X\rightarrow Y$ be any map.
A nonempty subset of $X$ of the form  $\varphi ^{-1}\left(\{y\}\right)$  will be called the {\it level set} (level line, level surface, level manifold...) of
$\varphi $ at level  $y\in Y$.

\medskip

{\tfp 2. The geometric inscribable polygons for $n\geq 4$}

\medskip

\noindent  First note that, following our definition 1.3, an inscribable  polygon $ (M_1, \cdots, M_n) $ is always convex. It is therefore
a star-shaped polygon with respect to any of its vertices. We can therefore represent it by an element
of ${\cal P}_n^\ast $. In this way, the set $\Gamma _n$ of inscribable polygons of ${\Bbb E}$ can be viewed
as a subset of the open set
$\Omega _n$. (Recall that any regular geometric  polygon is inscribable.)

\smallskip

\noindent {\tfpp 2.1. Remark.} The fact that a polygon $(M_1,\cdots ,M_n)$ is inscribable is equivalent to the fact that each one of the
$n-3$ quadrilaterals: $Q_1=(M_n,M_1,M_2,M_3), \;  \cdots ,
Q_{n-3}=(M_n,M_{n-3},M_{n-2},M_{n-1})$
is inscribable.
\smallskip

\noindent {\tfpp 2.2. Lemma.} {\it A convex quadrilateral $(A,B,C,D)$ is inscribable if and only if the distances $a=AB, \; b=BC, \; c=CD, \; d=DA\; $ and $\; e=BD\; $ satisfy the
relation}:
$$ad(b^2+c^2-e^2)+bc(a^2+d^2-e^2)=0.\leqno{(2.1)}$$

\smallskip

\noindent {\ib Proof.} Let $\alpha $ and $\beta $ denote the  measures respectively of the angles $\widehat{A}$   and   $\widehat{C}$. Then (by
the cosine's law of Al-Kashi):
$$0<\alpha <\pi ,\;  0<\beta  <\pi \hskip0.3cm
\hbox{and} \hskip0.3cm \displaystyle a^2+d^2-2ad\cos \alpha =e^2=b^2+c^2-2bc\cos
\beta .$$
We deduce:
$$\cases{\cos \alpha = {{a^2+d^2-e^2}\over {2ad}}\cr
{} \cr
\cos \beta ={{b^2+c^2-e^2}\over {2bc}}.} $$
On the other hand, the quadrilateral $
(A,B,C,D)$ is inscribable if and only if  $\alpha
+\beta  =\pi $ or, equivalently, if $\cos \beta
=-\cos \alpha $. Hence:
$$(A,B,C,D) \hbox{ inscribable} \Longleftrightarrow {{b^2+c^2-e^2}\over {2bc}}=-{{a^2+d^2-e^2}\over {2ad}}$$
which is also equivalent to $ad(b^2+c^2-e^2)+bc(a^2+d^2-e^2)=0.$

\smallskip

Remark 2.1 and Lemma 2.2 make possible to realize  $
\Gamma _n$ as a level set of a differentiable map. More precisely, consider the maps $\Theta :{\Bbb R}^5\longrightarrow {\Bbb R}$, $\gamma _k :\Omega _n\longrightarrow {\Bbb R}$
and
$\gamma  :\Omega _n\longrightarrow {\Bbb R}^{n-3}$ defined by:
$$\Theta (u)=u_1u_2(u_4^2+u_5^2-u_3^2)+u_4u_5(u_1^2+u_2^2-u_3^2)\hskip0.2cm \hbox{for} \hskip0.2cm
u=(u_1,u_2,u_3,u_4,u_5)\in {\Bbb R}^5,$$
$$\gamma _k (\omega )=\Theta (t_k,x_k,t_{k+1},x_{k+1},t_{k+2})
\hskip0.2cm \hbox{for} \hskip0.2cm k\in \{1,\cdots ,n-3\}$$ and:
$$\gamma  (\omega )=\left( \gamma _1 (\omega
),\cdots ,\gamma _{n-3} (\omega ) \right) \hskip0.2cm \hbox{where} \hskip0.2cm \omega = (t_1,x_1,t_2,\cdots ,x_{n-2},t_{n-1})\in \Omega_n.$$

\smallskip

\noindent {\tfpp 2.3. Proposition.} {\it   The set  $\Gamma _n$ of geometric inscribable  polygons is given by
$\Gamma _n=\gamma^{-1}\left(\{0\}\right)$. Moreover,  this set
is a differentiable submanifold of dimension $n$ of the Euclidean space  ${\Bbb R}^{2n-3}$.}

\smallskip

\noindent {\ib Proof.}   Let $\omega =(t_1,x_1,....,t_{n-2},x_{n-2},t_{n-1})$ be an element of $\Omega _n$  represented by a polygon $(M_1,\cdots ,M_n)$  (in ${\cal P}_n^\ast $). Then we have the following equivalences:

$\omega \in \Gamma_n  \Longleftrightarrow \hbox{$(M_n,M_k,M_{k+1},M_{k+2})$}\hbox{ is inscribable for $1\leq k\leq n-3$}$

\hskip1.2cm $\Longleftrightarrow \Theta
(t_k,x_k,t_{k+1},x_{k+1},t_{k+2})=0 \; \; \hbox{for $1\leq
k\leq n-3$}$

\hskip1.2cm $\Longleftrightarrow \gamma _k
(\omega )=0 \; \; \hbox{for $1\leq k\leq n-3$}$

\hskip1.2cm $\Longleftrightarrow \gamma  (\omega )=0$

\hskip1.2cm  $\Longleftrightarrow \omega \in
\gamma ^{-1}\left(\{0\}\right) .$

Now we will prove, by induction on  $n\geq 4$, that the map $\gamma $ has maximal rank at each point  $\omega $ of $\Omega _n$.

\noindent {\bf $\bullet \; $  The case $n=4$.}   We have $\gamma =\gamma _1$ and the map:
$$\gamma :\Omega _4\longrightarrow {\Bbb R}, \; \omega =(t_1,x_1,t_2,x_2,t_3)\longmapsto \gamma  (\omega )=t_1x_1(t_3^2+x_2^2-t_2^2)+t_3x_2(t_1^2+x_1^2-t_2^2)$$
has maximal rank at each point $\omega \in \Omega
_4$. Indeed, $\displaystyle {{\partial \gamma }\over {\partial
t_2}}(\omega )=-2t_2(t_1x_1+t_3x_2)\neq 0$. So $\Gamma
_4=\gamma ^{-1}\left(\{0\}\right)$  is a codimension $1$ submanifold of the open set
$\Omega _4$ and then a submanifold of dimension $4$ of
${\Bbb R}^5$.

\smallskip

\noindent $\bullet $ {\tfpp Heredity.}  Suppose that, for a fixed integer  $n\geq 4$, the map $\gamma $ has maximal rank at each
point of $\Omega _n$. If to each element $\omega
=(t_1,x_1,\cdots ,t_{n-2},x_{n-2},t_{n-1},x_{n-1},t_n)\in \Omega_{n+1}$ we associate:
$$\cases{
\omega' =(t_1,x_1,\cdots ,t_{n-2},x_{n-2},t_{n-1})\in \Omega_n\cr
{}\cr
\omega'' =(t_{n-2},x_{n-2},t_{n-1},x_{n-1},t_n)\in \Omega _4}$$
then we can write
$\omega=(\omega',x_{n-1},t_n)$ and $\gamma
(\omega)=\left(\gamma (\omega') , \gamma (\omega'')\right) .$
Note that here we are using the same notation  $
\gamma $ for three different maps. The Jacobian matrix of the map
$\gamma :\Omega _{n+1}\subset
{\Bbb R}^{2n-1}\rightarrow {\Bbb R}^{n-2}$  at each point $\omega \in
\Omega _{n+1}$  is given by:
$$\pmatrix{A(\omega')&0\cr B(\omega'')&C(\omega'')}\leqno{(2.2)}$$
where :
$$A(\omega')=\pmatrix{
{{\partial \gamma_1}\over {\partial t_1}}(\omega ')&{{\partial \gamma _1}\over {\partial x_1}}(\omega ')& \cdots &
{{\partial \gamma _1}\over {\partial t_{n-2}}}(\omega ')&{{\partial \gamma _1}\over {\partial x_{n-2}}}(\omega ')& {{\partial \gamma _1}\over {\partial t_{n-1}}}(\omega ')\cr
{}\cr
{{\partial \gamma _2}\over {\partial t_1}}(\omega ')&{{\partial \gamma _2}\over {\partial x_1}}(\omega ')& \cdots &
{{\partial \gamma _2}\over {\partial t_{n-2}}}(\omega ')&{{\partial \gamma _2}\over {\partial x_{n-2}}}(\omega ')& {{\partial \gamma _2}\over {\partial t_{n-1}}}(\omega ')\cr
{}\cr
\cdots & \cdots &\cdots &\cdots &\cdots &\cdots \cr {}\cr
{{\partial \gamma_{n-3}}\over {\partial t_1}}(\omega ')&{{\partial \gamma_{n-3}}\over {\partial x_1}}(\omega ')& \cdots &
{{\partial \gamma_{n-3}}\over {\partial t_{n-2}}}(\omega ')&{{\partial \gamma_{n-3}}\over {\partial x_{n-2}}}(\omega ')& {{\partial \gamma_{n-3}}\over {\partial t_{n-1}}}(\omega ')}
$$
which is the Jacobian matrix of the map
 $\gamma :\Omega _n\subset
{\Bbb R}^{2n-3}\rightarrow {\Bbb R}^{n-3}$  at the point $\omega '$,
$$B(\omega'')= \left( 0, \cdots , 0, {{\partial \Theta}\over {\partial u_1}}(\omega''), {{\partial \Theta}\over {\partial u_2}}(\omega''), {{\partial \Theta}\over {\partial u_3}}(\omega'')\right) \hskip0.2cm \hbox{and}\hskip0.2cm C(\omega'')=
\left( {{\partial \Theta}\over {\partial u_4}}(\omega''),{{\partial \Theta}\over {\partial u_5}}(\omega'')\right) .$$
By the induction hypothesis, the matrix $A(\omega ')$ has rank
$n-3$. On the other hand, the two partial derivatives:
$$\cases{
{{\partial \Theta }\over {\partial u_4}}(\omega '')
=2x_{n-1}t_{n-2}x_{n-2}+t_n\left(x_{n-2}^2+t_{n-2}^2-t_{n-1}^2\right)\cr
{}\cr
{{\partial \Theta }\over {\partial u_5}}(\omega
'')=2t_nt_{n-2}x_{n-2}+x_{n-1}\left(x_{n-2}^2+t_{n-2}^2-t_{n-1}^2\right)
} $$
do not vanish  simultaneously since otherwise we will have:
$$0=x_{n-1}{{\partial \Theta }\over {\partial u_4}}(\omega '')-t_n{{\partial \Theta }\over {\partial u_5}}(\omega '')=2x_{n-1}t_{n-2}(x_{n-1}^2-t_n^2)$$
which implies $x_{n-1}=t_n$. Taking into account this equality in the  relation  ${{\partial \Theta
}\over {\partial u_4}}(\omega '')=0$  we obtain  :
$$0=t_n\left[(x_{n-2}+t_{n-2})^2-t_{n-1}^2\right]$$
and then $x_{n-2}+t_{n-2}=t_{n-1}$. This contradicts the inequality $
t_{n-1}<x_{n-2}+t_{n-2}$ since these are the lengths of the sides of the non degenerate  triangle $(M_{n+1},M_{n-2},M_{n-1})$.

The Jacobian matrix of the map $\gamma :\Omega _{n+1}\subset {\Bbb R}^{2n-1}\rightarrow {\Bbb R}^{n-2}$ is then of rank $n-2$.
The proof by induction is then over.

We deduce that, for any $n\geq 4$, the map  $
\gamma :\Omega _n\subset {\Bbb R}^{2n-3}\rightarrow {\Bbb R}^{n-3}$ has maximal rank and then the nonempty set $\Gamma
_n=\gamma ^{-1}\left(\{0\}\right)$  is a codimension $(n-3)$ submanifold of the open set  $\Omega _n$ of $
{\Bbb R}^{2n-3}$. This implies that $\Gamma _n$ is a submanifold of  dimension $n$ of ${\Bbb R}^{2n-3}$.

\medskip

{\tfp 3. The area and perimeter foliations}
\medskip

\noindent In all this section the space ${\cal P}_n^\ast $ of star-shaped polygons will be identified
(as we have proceed until now) to the open set $\Omega_n$.

\smallskip

\noindent {\tfpp 3.1.} We define the maps $p : \Omega_n\longrightarrow {\Bbb R}$,
${\cal A} : \Omega_n\longrightarrow {\Bbb R}$ and $\Psi :\Omega_n\longrightarrow {\Bbb R}^2$ by:
$$\cases{p(\omega )=\hbox{perimeter of $\omega $}\cr {}\cr {\cal A}(\omega )= \hbox{area of $\omega $}\cr
{}\cr \Psi (\omega )=(p(\omega ),{\cal A}(\omega )).}\leqno{(3.1)}$$

$$\includegraphics[height=5.3cm]{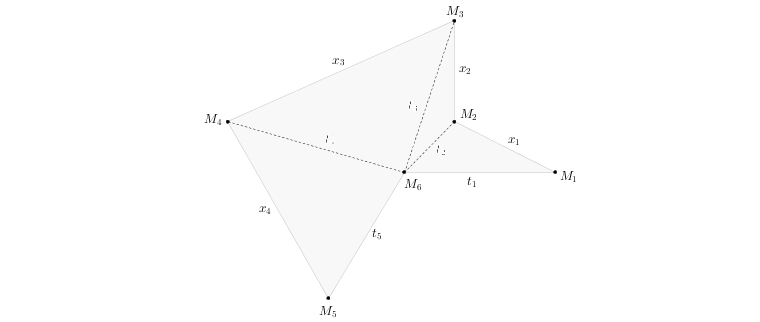}$$

\smallskip

For any $\omega =(t_1,x_1,t_2,\cdots ,t_{n-2},x_{n-2},t_{n-1})\in
\Omega _n$, we have:
$$\displaystyle \cases{
p(\omega )=t_1+x_1+x_2+\cdots +x_{n-2}+t_{n-1}\cr {}\cr
{\cal A}(\omega )={1\over 4}\sqrt{f(\omega_1)}+\cdots +{1\over 4}\sqrt{f(\omega _{n-2})}}\leqno{(3.2)}$$
with:
$$\cases{\omega _k=(t_k,x_k,t_{k+1})\in
{\cal V} \hbox{ for $k\in\{1,...,n-2\} $}\cr {}\cr
\hbox{area of $\omega_k ={1\over 4}\sqrt{f(\omega_k)}$ for $k\in \{ 1,\cdots ,n-2\}$ (H\'eron's formula)} \cr {}\cr
 f(x,y,z)=(x+y+z)(-x+y+z)(x-y+z)(x+y-z) \hbox{ for $(x,y,z)\in {\cal V}$}.}
$$
Setting  $s(x,y,z)=x+y+z$, we obtain for any $v=(x,y,z)\in {\cal V}$:
$$f(v)=s(v)\left(s(v)-2x\right)\left(s(v)-2y\right)\left(s(v)-2z\right) .$$

On the other hand, the maps $p, \; {\cal A}$ and $
\Psi $ are clearly  differentiable with gradient vectors
$\nabla p(\omega )$ and $\nabla \cal{A}(\omega )$
given by:
$$\cases{
\nabla p(\omega )=\left({{\partial p}\over {\partial t_1}}(\omega
),{{\partial p}\over {\partial x_1}}(\omega ),\cdots ,{{\partial
p}\over {\partial x_{n-2}}}(\omega ),
{{\partial p}\over {\partial t_{n-1}}}(\omega )\right)=(1,1,0,1,\cdots ,0,1,1)\cr {}\cr
\nabla {\cal A}(\omega )=\left( {{\nabla f(\omega
_1)}\over {8\sqrt{f(\omega _1)}}},\cdots ,{{\nabla f(\omega
_{n-2})}\over {8\sqrt{f(\omega _{n-2})}}}\right) =\left( {{\sqrt{f(\omega
_1)}}\over {8}}\cdot {{\nabla f(\omega _1)}\over {f(\omega
_1)}},\cdots ,{{\sqrt{f(\omega _{n-2})}}\over {8}}\cdot {{\nabla f(\omega
_{n-2})}\over {f(\omega _{n-2})}}\right)
}$$
where the logarithmic derivative ${{\nabla f}\over f}$
is given at each point $v=(x,y,z)\in {\cal V}$ by:
$${{\nabla f(v)}\over {f(v)}}={{\nabla s(v)}\over {s(v)}}+{{\nabla \left(s-2x\right)(v)}\over {s(v)-2x}}
+{{\nabla \left(s-2y\right)(v)}\over {s(v)-2y}}+{{\nabla
\left(s-2z\right)(v)}\over {s(v)-2z}}.$$
We then deduce:
$$\cases{
{1\over {f(v)}}{{\partial f}\over {\partial x}}(v)={1\over {s(v)}}-{1\over {s(v)-2x}}+{1\over {s(v)-2y}}+{1\over {s(v)-2z}}\cr {}\cr
{1\over {f(v)}} {{\partial f}\over {\partial y}}(v)={1\over {s(v)}}+{1\over {s(v)-2x}}-{1\over {s(v)-2y}}+{1\over {s(v)-2z}}\cr {}\cr
{1\over {f(v)}}{{\partial f}\over {\partial
z}}(v)={1\over {s(v)}}+{1\over {s(v)-2x}}+{1\over {s(v)-2y}}-{1\over {s(v)-2z}}.}$$
This gives the partial derivatives of ${\cal A}$:
$$\cases{
{{\partial {\cal A}}\over {\partial t_1}}(\omega )={{\sqrt{f(\omega _1)}}\over 8}\left( {1\over {s(\omega_1)}}-
{1\over {s(\omega _1)-2t_1}}+{1\over {s(\omega_1)-2x_1}}+{1\over {s(\omega_1)-2t_2}}\right) \cr
{}\cr
{{\partial {\cal A}}\over {\partial x_1}}(\omega )={{\sqrt{f(\omega _1)}}\over 8}\left({1\over {s(\omega_1)}}+{1\over {s(\omega_1)-2t_1}}-
{1\over {s(\omega_1)-2x_1}}+{1\over {s(\omega_1)-2t_2}}\right) \cr
{}\cr
{{\partial {\cal A}}\over {\partial t_{n-1}}}(\omega )={{\sqrt{f(\omega _{n-2})}}
\over 8}\left( {1\over {s(\omega_{n-2})}}+{1\over {s(\omega_{n-2})-2t_{n-2}}}+{1\over {s(\omega_{n-2})
-2x_{n-2}}}-{1\over {s(\omega_{n-2})-2t_{n-1}}}\right) }$$
and for $k\in \{ 2,\cdots ,n-2\} $:
$$\cases{ {{\partial {\cal A}}\over {\partial t_k}}(\omega )=
{{\sqrt{f(\omega _{k-1})}}\over 8}\left( {1\over {s(\omega_{k-1})}} +{1\over {s(\omega _{k-1})-2t_{k-1}}} +{1\over {s(\omega _{k-1})-2x_{k-1}}}
 -{1\over {s(\omega_{k-1})-2t_k}}
\right) \cr
{}\cr
\hskip 1.5cm +{{\sqrt{f(\omega _k)}}\over 8} \left( {1\over {s(\omega_k)}}+{1\over {s(\omega_k)-2x_k}}+{1\over {s(\omega _k)-2t_{k+1}}}-{1\over {s(\omega _k)-2t_k}}\right) \cr
{}\cr
{{\partial {\cal A}}\over {\partial x_k}}(\omega )={{\sqrt{f(\omega _k)}}\over 8}\left( {1\over {s(\omega_k)}}+{1\over {s(\omega _k)-2t_k}}
-{1\over {s(\omega_k)-2x_k}}+{1\over {s(\omega _k)-2t_{k+1}}} \right) .}$$

\noindent {\tfpp 3.2. Theorem.} {\it We have the following assertions.}

(1) {\it The perimeter function $p$ and the area function ${\cal A}$  are submersions on $
\Omega _n$. Then the level sets of
$p$ (resp. of ${\cal A}$) are leaves of a codimension $1$ foliation ${\cal F}_p$ (resp. ${\cal F}_a$) on $\Omega _n$.}

(2) {\it For $\omega \in \Omega _n$, the differential $d\Psi
(\omega )$ is of rank $2$ if $\omega $ is not a regular polygon and of rank $1$ if $\omega $ is a regular polygon. Then the map
$\Psi $ defines a codimension $2$ foliation ${\cal F}$ on the open set $\Omega _{n}$  of $
{\Bbb R}^{2n-3}$ which consists of non regular polygons $\omega
$ of $\Omega _n$.}

\smallskip

\noindent  {\ib Proof.} Let $\omega =(t_1,x_1,\cdots ,t_{n-2},x_{n-2},t_{n-1})$ be an element of $\Omega _n$  and
$(M_1,\cdots ,M_n)$ one of its representatives as a satr-shaped polygon.

\smallskip

\underbar{Point} (1)
\smallskip

\noindent {\bf $(\star )$} For any $\omega \in \Omega_n$, $dp(\omega )\neq 0$ since ${{\partial p}\over {\partial
t_1}}(\omega )=1\neq 0$. Then $p$  is a submersion on $\Omega _n$.

\smallskip

\noindent {\bf $(\star )$}  For any $\omega \in \Omega
_n$, ${{\partial {\cal A}}\over {\partial t_1}}(\omega )\neq 0$
or ${{\partial {\cal A}}\over {\partial x_1}}(\omega )\neq 0$. Indeed we have the  implications:
$$\cases{
{{\partial {\cal A}}\over {\partial t_1}}(\omega )=0\cr
\hbox{et}\cr
{{\partial {\cal A}}\over {\partial x_1}}(\omega )=0} \Longrightarrow  {{\partial {\cal A}}\over {\partial t_1}}(\omega )+{{\partial
{\cal A}}\over {\partial x_1}}(\omega )=0 \Longrightarrow
{2\over {s(\omega _1)}}+{2\over {s(\omega _1)-2t_2}}=0
\Longrightarrow t_1+x_1=0.$$
But the equality $t_1+x_1=0$ can not be satisfied.  Then $d{\cal A}(\omega )\neq 0$. This proves that ${\cal A}$  is a submersion on $\Omega _n$.
\smallskip

\underbar{Point} (2)
\smallskip

\noindent {\bf $(\star )$} If the sides of the polygon are not all equal, there exist two successive sides
having $M$ as common point and different lengths. One can suppose that $M=M_1$. In these conditions, the lengths $t_1=M_nM_1$ and $x_1=M_1M_2$ are
different. This implies
$\displaystyle {{\partial
{\cal A}}\over {\partial t_1}}(\omega )\neq {{\partial
{\cal A}}\over {\partial x_1}}(\omega )$   since we have the
implication ${{\partial {\cal A}}\over {\partial t_1}}(\omega )={{\partial {\cal A}}\over {\partial x_1}}(\omega )\Longrightarrow t_1=x_1.$
The Jacobian matrix of the map $\Psi $ at the point $\omega $ is:
$${\cal J}(\Psi ,\omega )=\pmatrix{1&1&0&1&\cdots &0&1&1\cr {}\cr
{{\partial {\cal A}}\over {\partial t_1}}(\omega )&{{\partial {\cal A}}\over {\partial x_1}}(\omega )&{{\partial {\cal A}}\over {\partial t_2}}(\omega )&
{{\partial {\cal A}}\over {\partial x_2}}(\omega )&\cdots &{{\partial {\cal A}}\over {\partial t{n_2}}}(\omega )&{{\partial {\cal A}}\over {\partial x_{n-2}}}(\omega )&
{{\partial {\cal A}}\over {\partial t_{n-1}}}(\omega )} $$
and then it is of rank $2$ since the $2\times 2$-matrix of  ${\cal J}(\Psi ,\omega )$ consisting of the two first columns
is invertible.

\smallskip

\noindent {\bf $(\star )$} Suppose the polygon is equilateral. We consider two cases:

\smallskip

$\bullet $ The condition $(C)$ below is satisfied.

\smallskip
$$\cases{\hbox{The partial derivatives }
{{\partial {\cal A}}\over {\partial t_1}}(\omega ), {{\partial {\cal A}}\over {\partial t_{n-1}}}(\omega ), {{\partial {\cal A}}\over {\partial x_k}}(\omega )\cr
\hbox{ are not equal, $k\in \{ 1,\cdots ,n-2\} $}\cr
\hbox{or}\cr
\hbox{at least one of the derivatives ${{\partial {\cal A}}\over {\partial t_k}}(\omega )$, $k\in \{ 1,\cdots ,n-2\} $ is not zero.} }\leqno{\hbox{(C)}}$$
The matrix ${\cal J}(\Psi ,\omega )$ has rank
$2$ since it admits a matrix of order $2$ which is reversible.

\smallskip

$\bullet $ The  condition $(C)$ is not satisfied.
\smallskip
In this case the following  condition non(C) is satisfied:

$$\cases{\hbox{The partial derivatives }
{{\partial {\cal A}}\over {\partial t_1}}(\omega ), {{\partial {\cal A}}\over {\partial t_{n-1}}}(\omega ), {{\partial {\cal A}}\over {\partial x_k}}(\omega )\cr
\hbox{ are all equal for $k\in \{ 1,\cdots ,n-2\} $}\cr
\hbox{and}\cr
{{\partial {\cal A}}\over {\partial t_k}}(\omega ), \hbox{ $k\in \{ 1,\cdots ,n-2\} $ are all zero} }\leqno{\hbox{non(C)}}$$
and we have:
$${\cal J}(\Psi ,\omega )=\pmatrix{1&1&0&1&\cdots &0&1&1\cr {}\cr
{{\partial {\cal A}}\over {\partial t_1}}(\omega )&{{\partial {\cal A}}\over {\partial t_1}}(\omega )&0&
{{\partial {\cal A}}\over {\partial t_1}}(\omega )&\cdots &0&{{\partial {\cal A}}\over {\partial t_1}}(\omega )&
{{\partial {\cal A}}\over {\partial t_1}}(\omega )} \leqno{(3.3)}$$
This implies that ${\cal J}(\Psi ,\omega )$ is of rank $1$. We shall prove that, in this case,
$\omega $ is regular polygon.

\smallskip

For $k\in \{2,\cdots , n-2\}$, we have  $\displaystyle
{{\partial {\cal A}}\over {\partial x_{k-1}}}(\omega
)={{\partial {\cal A}}\over {\partial x_k}}(\omega )$ which implies:
$$\displaystyle \left({{\partial{\cal A}}\over {\partial
x_{k-1}}}(\omega )\right)^2-\left({{\partial
{\cal A}}\over {\partial x_k}}(\omega )\right)^2=0.$$
Thus, taking into account the fact that the sides  $x_k$ are all equal in this case, we
obtain by  factorization :
$${{64 x_k^2 t_k^2 (t_{k-1} - t_{k+1})(t_{k-1} +
t_{k+1})(t_{k-1}t_{k+1} - x_k^2 + t_k^2)(t_{k-1} t_{k+1} + x_k^2 -
t_k^2)}\over \Delta }=0$$
where:
$$\begin{array}{rcl}
\Delta = &(x_k - t_k - t_{k+1})(x_k - t_k + t_{k+1})(x_k + t_k -
t_{k+1})(x_k + t_k + t_{k+1})\cr &(t_{k-1} - x_k - t_k) (t_{k-1} - x_k +
t_k)(t_{k-1} + x_k - t_k)(t_{k-1} + x_k + t_k).
\end{array}$$
Thus $\displaystyle (t_{k-1} -
t_{k+1})(t_{k-1}t_{k+1} - x_k^2 + t_k^2)(t_{k-1} t_{k+1} + x_k^2 -
t_k^2)=0.$

\smallskip

\noindent {\bf $(\star )$}  If $t_{k-1}-t_{k+1}=0$, taking into account  the
relation ${{\partial {\cal A}}\over {\partial t_k}}(\omega
)=0$, we obtain:
$$t_k^2=x_k^2+t_{k+1}^2.$$
This implies that the two  triangles $(M_n,M_{k-1},M_k)$ and $(M_n,M_{k+1},M_k)$ are rectangular respectively at $M_{k-1}$ and $M_{k+1}$.
Consequently  the  quadrilateral $(M_n,M_{k-1},M_k,M_{k+1})$ is inscribable ($\widehat{M}_{k-1}={\pi \over 2}=\widehat{M}_{k+1}$).
\smallskip

\noindent {\bf $(\star )$} If $t_{k-1}t_{k+1} - x_k^2 + t_k^2=0$, the convex quadrilateral  $(M_n,M_{k-1},M_k,M_{k+1})$ has the following properties:
$$\cases{
M_{k-1}M_k=x_k=M_kM_{k+1}\; \hbox{ (The polygon is equilateral in this case.)}\cr
{}\cr
t_k^2=x_k^2-t_{k-1}t_{k+1}\cr
{}\cr
\cos
\widehat{M_{k-1}}={{t_k^2-t_{k-1}^2-x_k^2}\over {2t_{k-1}x_k}}={{x_k^2-t_{k-1}t_{k+1}-t_{k-1}^2-x_k^2}\over {2t_{k-1}x_k}}=-{{t_{k+1}+t_{k-1}}\over {2x_k}}\cr
{}\cr
\cos
\widehat{M_{k+1}}={{t_k^2-t_{k+1}^2-x_k^2}\over {2t_{k+1}x_k}}{{x_k^2-t_{k-1}t_{k+1}-t_{k+1}^2-x_k^2}\over {2t_{k+1}x_k}}=-{{t_{k-1}+t_{k+1}}\over {2x_k}}\cr
{}\cr
\widehat{M_{k-1}}=\widehat{M_{k+1}} \hbox{ (by the equality of the cosines).}
}\leqno{(3.4)}$$
The triangle $(M_{k-1},M_n,M_{k+1})$ is then isosceles
at the $M_n$ and $t_{k-1}=t_{k+1}$. This implies, like in the preceding case, that the quadrilateral $
(M_n,M_{k-1},M_k,M_{k+1})$ is inscribable.

\smallskip

\noindent {\bf $(\star )$} If $\displaystyle t_{k-1} t_{k+1} + x_k^2 -
t_k^2=0$, then the convex quadrilateral $(M_n,M_{k-1},M_k,M_{k+1})$ is still  inscribable since it satisfies the
relation  $\cos \widehat{M_{k-1}}=-\cos
\widehat{M_{k+1}}$. Indeed, we have:
$$\cases{
t_k^2=t_{k-1}t_{k+1} + x_k^2\cr
{}\cr
\cos
\widehat{M_{k-1}}={{t_k^2-t_{k-1}^2-x_k^2}\over {2t_{k-1}x_k}}={{x_k^2+t_{k-1}t_{k+1}-t_{k-1}^2-x_k^2}\over {2t_{k-1}x_k}}={{t_{k+1}-t_{k-1}}{2x_k}}\cr
{} \cr
\cos
\widehat{M_{k+1}}={{t_k^2-t_{k+1}^2-x_k^2}\over {2t_{k+1}x_k}}{{x_k^2+t_{k-1}t_{k+1}-t_{k+1}^2-x_k^2}\over {2t_{k+1}x_k}}={{t_{k-1}-t_{k+1}}\over {2x_k}}.}\leqno{(3.5)}
$$

We have proved that, in all cases, the quadrilateral
$(M_n,M_{k-1},M_k,M_{k+1})$ is inscribable for any  $k\in \{2,\cdots ,n-2\}$.
Then the polygon    $(M_1,\cdots ,M_n)$ is inscribable. Since the latter is equilateral, it is necessarily regular.

\smallskip

Finally the singular points of the map  $\Psi $ are the regular polygons and
this map induces a submersion on  $\Omega _n^\ast $ whose level sets are leaves of a foliation  ${\cal F}$.

\medskip

{\tfp 4. The example of triangles}

\medskip

\noindent It is the situation where we see things more concretely and where the drawings are more visible. The way to treat the topic in this section
will be slightly different from the others.

\smallskip

\noindent {\tfpp 4.1. The space of non degenerate triangles}
\smallskip

To give oneself a non degenerate triangle (in any Euclidean finite dimensional space) is to give oneself three real numbers $x>0$, $y>0$ and $z>0$ such that:
$$\cases{x<y+z\cr y<z+x \cr z<x+y}\leqno{(4.1)}$$
which represent the lengths of the sides. Exceptionally in this section, we shall denote a triangle by
$\langle xyz\rangle $ instead of
$(X,Y,Z)$ where the points $X$, $Y$  and $Z$ are the vertices. Indeed it is well known that  $\langle xyz\rangle $ is isometric to $\langle x'y'z'\rangle $ if
$x=x'$, $y=y'$ and $z=z'$. (For the moment we will make the difference between a triangle and another obtained by permutation of the three numbers representing it
even if, geometrically, they are the same!)
From now on,  $\lambda $ will be the half perimeter $\lambda ={{x+y+z}\over 2}$.
\smallskip

The set of non degenerate triangles  is thus the open set $\Omega_3\subset ({\Bbb R}_+^\ast )^3$ given by (1.7).
We will describe it explicitly. To inequalities (4.1) are associated three equations
defining respectively three planes:
$$\cases{\Sigma_1=\{ x=y+z\} \cr \Sigma_2=\{ y=z+x\} \cr \Sigma_3=\{ z=x+y\} } .\leqno{(4.2)}$$
In the slice $\{ x+y+z=2\lambda \} $ of ${\Bbb R}_+^3$, $\Sigma_1$, $\Sigma_2$ and $\Sigma_3$ are the sides of an equilateral triangle in ${\Bbb R}_+^3$
whose vertices are $X_\lambda =(0,\lambda ,\lambda )$, $Y_\lambda =(\lambda ,0,\lambda )$ and $Z_\lambda =(\lambda ,\lambda ,0)$ (see the picture bellow);
the interior $P_\lambda $ of the convex hull   of these three points
 represents the space
of triangles $\langle xyz\rangle $ whose perimeter is $2\lambda $.
\smallskip

When $\lambda $ varies to $\lambda'$, we obtain another  $P_{\lambda'}$, image of $P_{\lambda }$ by the homothety centered at the
origin and with ratio $k={{\lambda'}\over \lambda }$. Thus, the space  $\Omega_3$ is {\it foliated} by these $P_\lambda $; $\Omega_3$ is in fact the open cone
with vertex the origin and basis anyone of these plaques $P_{\lambda }$, for instance $P_1$:
$$\Omega_3=\bigcup_{\lambda \in {\Bbb R}_+^\ast }\lambda P_1=\left\{ \lambda X:\hbox{ $X\in P_1$  et $\lambda \in {\Bbb R}_+^\ast $}\right\} .\leqno{(4.3)}$$

$$\includegraphics[height=8cm]{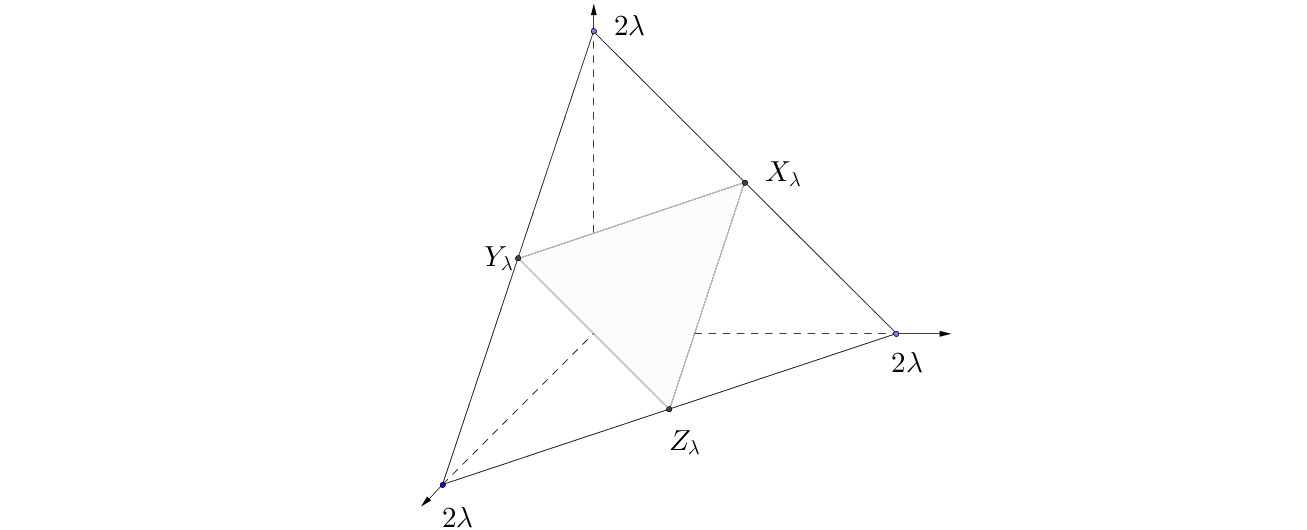}$$

For a particular situation which will appear thereafter, we recall the following result which we have already established in the general case of polygons.
\smallskip

\noindent {\it For a given family of  triangles  with prescribed perimeter, the maximum of the area is realized by the equilateral triangle.}

\smallskip

\noindent {\tfpp 4.2. The perimeter foliation ${\cal F}_p$}
\smallskip

Each $P_\lambda $ (where $\lambda \in {\Bbb R}_+^\ast $) is the level set $p(x,y,z)=2\lambda $ where $p$ is the perimeter function
$p(x,y,z)=x+y+z$. We have also seen that the level surface $P_\lambda $ is the interior of the convex hull of the triangle $X_\lambda Y_\lambda Z_\lambda $.
\smallskip

Thus we have a foliation ${\cal F}_p$ on $\Omega_3$ whose leaves are the surfaces $P_\lambda $ ($\lambda >0$).
Of course, ${\cal F}_p$ is trivial since isomorphic to the product
$P_1\times {\Bbb R}_+^\ast $.

\medskip
\noindent {\tfpp 4.3. The area foliation ${\cal F}_a$}
\smallskip

The function ${\cal A} :({\Bbb R}_+^\ast )^3\longrightarrow {\Bbb R}_+^\ast $ which associates to a triangle $\langle xyz\rangle $ its area is given by  H\'eron formula:
$${\cal A}(x,y,z)={1\over 4}\sqrt{(x+y+z)(-x+y+z)(x-y+z)(x+y-z)}.\leqno{(4.4)}$$
The foliation ${\cal F}_a$ by which we will be interested  is the foliation whose leaves are
the level surfaces of this function.

\smallskip

$\bullet $ The surface at level $s$ of the function ${\cal A}$ on the open set $\Omega_3$ is exactly the surface at level $s^2$ of the function
$\Phi ={\cal A}^2$. The benefit of working
with $\Phi $ instead of ${\cal A}$ is that there is no more square root, which simplifies the calculations, among others that of the differential which plays a fundamental role.
We consider then the function:
$$\Phi (x,y,z)={1\over 16}(x+y+z)(-x+y+z)(x-y+z)(x+y-z).\leqno{(4.5)}$$

$\bullet $ The differential of $\Phi $ has the form:
$$d\Phi (x,y,z)={1\over 16}\{ A(x,y,z)dx+B(x,y,z)dy+C(x,y,z)dz\} $$
where the functions $A$, $B$ and $C$
are given as follows:
$$\begin{array}{rcl}
A&= \hbox{$(-x+y+z)(x-y+z)(x+y-z) -(x+y+z) (x-y+z)(x+y-z)$}\cr
&\hskip0.4cm \hbox{$+(x+y+z)(-x+y+z)(x+y-z)+(x+y+z)(-x+y+z)(x-y+z)$} \cr
{}\cr
B&= \hbox{$(-x+y+z)(x-y+z)(x+y-z)  +(x+y+z) (x-y+z)(x+y-z)$}\cr
& \hskip0.4cm  \hbox{$-(x+y+z)(-x+y+z)(x+y-z)+(x+y+z)(-x+y+z)(x-y+z)$} \cr
{}\cr
C&= \hbox{$ (-x+y+z)(x-y+z)(x+y-z)+(x+y+z) (x-y+z)(x+y-z)$}\cr
& \hskip0.4cm \hbox{$+(x+y+z)(-x+y+z)(x+y-z)-(x+y+z)(-x+y+z)(x-y+z)$}.
\end{array} \leqno{(4.6)}$$

An easy but long computation shows that these three functions $A$, $B$ et $C$  are zero simultaneously only if $x=y=z=0$,
which can not happen since $ (0,0,0) $ is not in $\Omega_3$.

\smallskip

$\bullet $ If we fix the perimeter $2\lambda $, the area function $a$ is maximal, and so is the function $\Phi $, when
 $x=y=z={2\over 3}\lambda $;
at this point $\Phi $ is equal to  ${{\lambda^4}\over {27}}$. These are the values taken by the function $\Phi $ on the open half line $\Delta $ whose equation is $x=y=z$.
\smallskip

Now let  $\Omega_3^\ast $ be the open set $\Omega_3\setminus \Delta $. At $u=(x,y,z)\in \Omega_3^\ast $, the differential $d_u\Phi $ has rank   $1$; then
the set level of $\Phi $ passing through this point is a regular surface $A$, in fact an algebraic surface of degree $4$. Its equation is:
$$(x+y+z)(-x+y+z)(x-y+z)(x+y-z)=16\Phi (u).$$

\smallskip

Let $G$  be the subgroup of Isom$({\Bbb R}^3)$ (the full group of isometries of the Euclidean space ${\Bbb R}^3$) generated by the
rotation whose axis is  $\Delta $ and angle ${{2\pi }\over 3}$ and the reflection $\sigma $ with respect to the plane of equation $x=y$.
(The restrictions of these elements to the plane of equation $x+y+z=2\lambda $
is the group of isometries of the equilateral triangle
$X_\lambda Y_\lambda Z_\lambda $.)  It leaves the space $\Omega_3$ invariant and also its boundary $\partial \Omega_3$, the half line
$\Delta $ and the open sets $\Omega_3$ and $\Omega_3^\ast$. Then it acts on $\Omega_3$ and fixes each leaf
of ${\cal F}_a$; the same applies to the foliation ${\cal F}_p$.

\smallskip

\noindent {\tfpp 4.4.} Let $\Psi :\Omega_3^\ast \longmapsto \left( {\Bbb R}_+^\ast \right)^2$ be the function :
$$\Psi (x,y,z)=(p(x,y,z),\Phi (x,y,z)).$$
Up to a multiplicative factor, the matrix of its differential at  $u=(x,y,z)$ is :
$$d_u\Psi =\pmatrix{1&1&1\cr A(u)&B(u)&C(u)}$$
where $A$, $B$ and $C$ are the functions given by  (5.6).
It can be shown that these functions are equal only if $x=y=z$; then, for $u\in \Omega_3^\ast $,
$d_u\Psi $ has rank  $2$. Thus, the level sets of $\Psi $ are regular curves, leaves of a foliation
${\cal F}_\ast $ on $\Omega_3^\ast $.
\smallskip

\noindent {\tfpp 4.5.} On $\Omega_3$ we have a singular foliation   ${\cal F}={\cal F}_p\cap {\cal F}_a$. Its leaves of dimension $0$ are the
points of the open half line $\{ \left( {2\over 3}\lambda ,{2\over 3}\lambda ,{2\over 3}\lambda \right) :\lambda \in {\Bbb R}_+\} $. The other leaves
are of dimension  $1$; each one
has equation $\Psi (u)=\hbox{constant}$ in the open set $\Omega_3^\ast $. These curves define (by restriction) a foliation on each
plaque $P_\lambda $ (leaf of ${\cal F}_p$). To see what it is, this plaque is projected orthogonally on the plane $z = 0$; we obtain the
foliation drawn  on the picture below. We will explain what all this means.

\vskip0.4cm
$$\includegraphics[height=5cm]{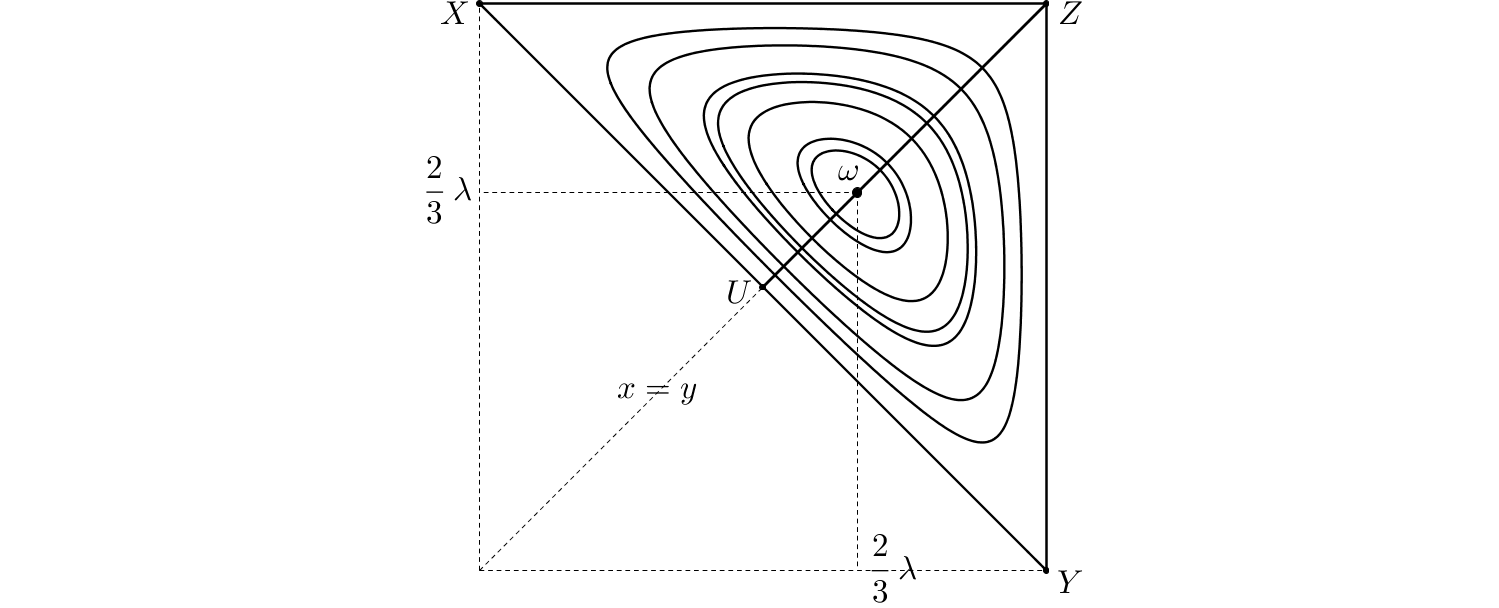}$$

\medskip

\noindent The interior of the triangle $XYZ$  is the projection (which we denote by $\Theta_\lambda $) on the plane $z=0$
of the set  $P_\lambda $ of triangles $\langle x_\lambda y_\lambda z_\lambda \rangle $ with perimeter $2\lambda $. Note that the boundary of
$P_\lambda $ is an equilateral triangle  while $XYZ$ is an isosceles and rectangle  triangle.
The foliation ${\cal F}$ on $P_\lambda $ is isomorphic to the foliation on the picture via the diffeomorphism $f:P_\lambda \longrightarrow \Theta_\lambda $
defined by $f(x,y,z)=(x,y,0)$ with inverse $f^{-1}(x,y,0)=(x,y,2\lambda -x-y)$.
\smallskip

$\bullet $ The point  $\omega $ with coordinates $\left( {2\over 3}\lambda ,{2\over 3}\lambda \right) $ corresponds to the equilateral triangle
$\langle xxx\rangle $ with maximal area. As we easily imagine, an equilateral triangle may  never be deformed to an other one  having the same perimeter and
the same area.
\smallskip

$\bullet $  The curves at the interior of $\Theta_\lambda $
are leaves of a foliation of $\Theta_\lambda  \setminus \{ \omega \} $, each leaf corresponds to the set of triangles having the same area.  It has
$\lambda (2\lambda -x)(2\lambda -y)(x+y)=8c$ as equation where $c$ is a constant varying in the interval  $\big\rbrack 0,{{8\lambda^4}\over {27}}\big\lbrack $.
\smallskip

$\bullet $ The piece $UZ$ of diagonal corresponds to isosceles triangles (for which $x=y$). In each leaf, there is exactly the projections of two isosceles triangles
$\langle xxz\rangle $ and $\langle x'x'z'\rangle $.
\smallskip

\noindent {\tfpp 4.6.} Geoffrey Letellier constructed two lines of isosceles triangles : $\langle x_\lambda y_\lambda z_\lambda \rangle $ and $\langle x_\lambda 'y_\lambda 'z_\lambda '\rangle $
where $\lambda \in {\Bbb R}_+^\ast $ with $x_\lambda =y_\lambda ={11\over 14}\lambda $, $z_\lambda ={3\over 7}\lambda $ and $x_\lambda '=y_\lambda '={4\over 7}\lambda $, $z_\lambda '={6\over 7}\lambda $. They are such that, for any $\lambda \in {\Bbb R}_+^\ast $:
$$\cases{\hbox{$\langle x_\lambda x_\lambda z_\lambda \rangle $ and $\langle x_\lambda 'x_\lambda 'z_\lambda '\rangle $
have the same perimeter  $2\lambda $.}\cr
{}\cr
\hbox{$\langle x_\lambda x_\lambda z_\lambda \rangle $ and $\langle x_\lambda 'x_\lambda 'z_\lambda '\rangle $ have the same area ${{3\lambda^2}\over {7\sqrt{7}}}$.}\cr
{}\cr
\hbox{$\langle x_\lambda x_\lambda z_\lambda \rangle $ and $\langle x_\lambda 'x_\lambda 'z_\lambda '\rangle $ are not isometric.}}$$

For instance, the two isosceles triangles $x=y=11$, $z=6$ and $x'=y'=8$, $z'=12$ have the same
perimeter equal to  $28$ and the same area equal to $12\sqrt{7}$.
\smallskip

$\bullet $ Finally one can see on the picture that all the situation is invariant by the reflection  $\sigma $
(symmetry with respect to the diagonal  $x=y$) while that on the triangle $P_\lambda $ is invariant by the full group $G$. \hfill $\diamondsuit $

\medskip

{\tfp 5.  Some results related to the perimeter and the area}

\medskip

\noindent The following well known classical results are among the most beautiful theorems that we can cite in Euclidean elementary geometry of the plane.

\smallskip

\noindent {\tfpp 5.1. Theorem (Isoperimetric inequality).} {\it Among all the convex polygons  with prescribed perimeter, the regular polygon is the one whose area
is maximum.}

\smallskip

For a sketch of proof, see for instance [Han]. In the same order of ideas, we also have the following theorem.
Its proof is not difficult but it  is a bit long and not immediate. (And the reader
can even attempt to reproduce it himself!)

\smallskip

\noindent {\tfpp 5.2. Theorem.} {\it Among all convex polygons
whose sides have given lengths, the inscribable polygon
is the one whose area is maximum.}

\smallskip

Using the analytic expression of the function ``area"
 ${\cal A}:\omega \in \Omega _n \longrightarrow \hbox{area}(\omega )\in {\Bbb R}$,
we prove the following result related to the two theorems
above. (It was also partially established, by a different method, in [Khi].)

\smallskip

\noindent {\tfpp 5.3. Theorem.} {\it We have the following results.}
\smallskip

(1) {\it For any real number $L>0$, the differentiable manifold  $p^{-1}\left(\{L\}\right)$ consisting of all polygons with perimeter  $L$,  is diffeomorphic
to a convex open set
 $\Omega_{n,L}$ of ${\Bbb R}^{2n-4}$ and the
restriction  ${\cal A}_L : \Omega _{n,L} \rightarrow {\Bbb R}$
of ${\cal A}$ to $\Omega _{n,L}$ admits a critical  point
at the unique  regular polygon   $\omega_L$ of
perimeter $L$.}
\smallskip

(2) {\it  The convex polygons whose sides have
given lengths form a differentiable manifold diffeomorphic
to an open convex set of   ${\Bbb R}^{n-3}$  and the restriction of the function ${\cal A}$ to this open set
 admits a critical point at its unique inscribable polygon.}
\smallskip

\noindent {\ib Proof.}  Recall that, for any $\omega
=(t_1,x_1,\cdots ,t_{n-2},x_{n-2},t_{n-1})\in \Omega _n$, we have:
$$ p(\omega )=\hbox{perimeter}(\omega )=t_1+x_1+\cdots +x_{n-2}+t_{n-1}$$
$${\cal A}(\omega )=\hbox{area}(\omega )={1\over 4}\sqrt{f(\omega _1)}+\cdots +{1\over 4}\sqrt{f(\omega
_{n-2})}$$
with:
$$\cases{\omega _k=(t_k,x_k,t_{k+1})\in {\cal V} \hbox{ for $k\in\{1,\cdots ,n-2\}$}\cr
{}\cr
 f(x,y,z)=(x+y+z)(-x+y+z)(x-y+z)(x+y-z) \hbox{ for $(x,y,z)\in {\cal V}$}. }
$$
Setting  $h(v)={1\over 4}\sqrt{f(v)}$ for $v\in {\cal V}$, we obtain:
$${\cal A}(\omega )=h(\omega _1)+\cdots +h(\omega _{n-2})$$
for any $\omega
=(t_1,x_1,\cdots ,t_{n-2},x_{n-2},t_{n-1})\in \Omega _n$.
\smallskip

\underbar{Point} (1)

\smallskip
Let $L\in ]0,+\infty [$. For any $\omega =(t_1,x_1,t_2,...,t_{n-2},x_{n-2},t_{n-1})\in \Omega _n$, we have:
$$ p(\omega )=L\Leftrightarrow t_1=L-x_1-\cdots -x_{n-2}-t_{n-1}$$
Then, by considering the affine map $T:{\Bbb R}^{2n-4}\rightarrow {\Bbb R}^{2n-3}={\Bbb R}\times
{\Bbb R}^{2n-4}$ given by:
$$u=(x_1,t_2,\cdots ,t_{n-2},x_{n-2},t_{n-1})\longmapsto
T(u)=(t_1(u),u)$$
where $t_1(u)=L-x_1-\cdots -x_{n-2}-t_{n-1}$,  we see that
$p^{-1}\left(\{L\}\right)$ is naturally identified to the convex open set  $\Omega _{n,L}=T^{-1}\left(\Omega _n\right)$ of ${\Bbb R}^{2n-4}$.
\smallskip

For any  $u=(x_1,t_2,\cdots ,t_{n-2},x_{n-2},t_{n-1})\in \Omega _{n,L}$, we have:
$${\cal A}_L(u)
=h(t_1(u),x_1,t_2)+h(t_2,x_2,t_3)+\cdots +h(t_{n-2},x_{n-2},t_{n-1}).$$
Set $\displaystyle t_1=t_1(u)$ and $\displaystyle \omega _k=(t_k,x_k,t_{k+1})$ for $k\in
\{1,\cdots ,n-2\}$. The partial derivatives for any $(x,y,z)\in {\cal V}$ are:
$$ \cases{
{{\partial h}\over {\partial x}}(x,y,z)={{x
\left(-x^2+y^2+z^2\right)}\over {2 \sqrt{(x+y+z))(-x+y+z)(x-y+z)
(x+y-z)}}}\cr
{}\cr
{{\partial h}\over {\partial y}}(x,y,z)={{y
\left(x^2-y^2+z^2\right)}\over {2 \sqrt{(x+y+z))(-x+y+z)(x-y+z)
(x+y-z)}}}\cr
{}\cr
{{\partial h}\over {\partial z}}(x,y,z)={{z\left(x^2+y^2-z^2\right)}\over {2 \sqrt{(x+y+z))(-x+y+z)(x-y+z)
(x+y-z)}}}
} \leqno{(5.1)}$$
$$\cases{
{{\partial {\cal A}_L}\over {\partial x_1}}(u)={{\partial
}\over {\partial x_1}}\left[h\left(t_1(u),x_1,t_2\right)\right]=-
{{\partial h}\over {\partial x}}(\omega _1)+{{\partial
h}\over {\partial y}}(\omega _1)=\cr
{}\cr
{{(t_1-x_1)(t_1+x_1-t_2)(t_1+x_1+t_2)}\over {2 \sqrt{(t_1+x_1+t_2)(-t_1+x_1+t_2)(t_1-x_1+t_2)(t_1+x_1-t_2)}}}\cr
{}\cr
{{\partial {\cal A}_L}\over {\partial t_{n-1}}}(u)={{\partial
}\over {\partial
t_{n-1}}}\left[h\left(t_1(u),x_1,t_2\right)+h\left(t_{n-2},x_{n-2},t_{n-1}\right)\right]=
-{{\partial h}\over {\partial x}}(\omega _1)+{{\partial
h}\over {\partial z}}(\omega _{n-2})}$$
and for $k\in \{2,\cdots ,n-2\}$:
$$\cases{
{{\partial {\cal A}_L}\over {\partial t_k}}(u)={{\partial
}\over {\partial
t_k}}\left[h\left(t_{k-1},x_{k-1},t_k\right)+h\left(t_k,x_k,t_{k+1}\right)\right]=
{{\partial h}\over {\partial z}}(\omega _{k-1})+{{\partial
h}\over {\partial x}}(\omega _k)\cr
{}\cr
{{\partial {\cal A}_L}\over {\partial x_k}}(u)={{\partial
}\over {\partial
x_k}}\left[h\left(t_1(u),x_1,t_2\right)+h\left(t_k,x_k,t_{k+1}\right)\right]=
-{{\partial h}\over {\partial x}}(\omega _1)+{{\partial
h}\over {\partial y}}(\omega _k)
}\leqno{(5.2)}$$

\smallskip

\noindent {\bf $(\star )$} If the sides of the polygon are not all
equal, there are two consecutive ones with a
common vertex $M$ and different lengths.

By changing the numbering of the vertices, one can
assume $M=M_1$. In these conditions, the lengths $
t_1=M_nM_1$ and $x_1=M_1M_2$ are different and this implies $\displaystyle {{\partial
{\cal A}_L}\over {\partial x_1}}(u)\neq 0$ and that $u$ is not a critical point of  ${\cal A}_L$.

\smallskip

\noindent {\bf $(\star )$}  If the lengths of all the sides are equal, then $t_1=x_1=x_2=...=x_{n-2}=t_{n-1}$   and
a necessary condition for this polygon to be a critical point of ${\cal A}_L$,
is ${{\partial {\cal A}_L}\over {\partial
t_k}}(u)=0$ for any $k\in \{2,\cdots ,n-2\} $, or:
$${{\partial h}\over {\partial
z}}(t_{k-1},x_k,t_k)=-{{\partial h}\over {\partial
x}}(t_k,x_k,t_{k+1}) \; \; \hbox{pour tout $k\in
\{2,\cdots ,n-2\} $}.$$
This implies:
$$\left({{\partial h}\over {\partial
z}}(t_{k-1},x_k,t_k)\right)^2-\left({{\partial h}\over {\partial
x}}(t_k,x_k,t_{k+1})\right)^2=0\; \; \hbox{for any $k\in
\{2,\cdots ,n-2\}$}.$$
The development of this relationship leads to next equality:
$${{(t_{k+1}t_{k-1} + x_k^2 - t_k^2)(t_{k+1}t_{k-1} - x_k^2 +
t_k^2)(t_{k+1} + t_{k-1})(t_{k+1} - t_{k-1})x_k^2t_k^2}\over \Sigma }=0$$
where:
$$\begin{array}{rcl}
\Sigma = &(t_{k+1} +
x_k + t_k)(t_{k+1} + x_k - t_k)(t_{k+1} - x_k + t_k)(t_{k+1} - x_k
- t_k)\cr &(t_{k-1} + x_k + t_k)(t_{k-1} + x_k - t_k)(t_{k-1} - x_k +t_k)(t_{k-1} - x_k - t_k).
\end{array}$$
Thus $\displaystyle (t_{k+1} -
t_{k-1})(t_{k+1}t_{k-1} - x_k^2 + t_k^2)(t_{k+1}t_{k-1} + x_k^2 -
t_k^2)=0.$
The proof ends as that of the Theorem 3.2.
We thus obtain the inscriptibility of all the quadrilaterals
$(M_n,M_{k-1},M_k,M_{k+1})$ for $k\in \{2,\cdots ,n-2\}$ and then the inscriptibility of the polygon   $(M_1,\cdots ,M_n)$.
But since the latter has all its sides of the same length, it is necessarily regular.

\smallskip

Finally the singular points of the map   ${\cal A}_L$ are the regular polygons of   $\Omega_{n,L}$,  that is, the unique regular
polygon  $\omega_L$ of perimeter  $L$.

\smallskip

\underbar{Point} (2)

\smallskip

Let  $v
=(\overline{t}_1,\overline{x}_1,\cdots ,\overline{x}_{n-2},\overline{t}_{n-1})\in
{\cal V}_n$ and let $F_v$ be the set of convex polygons whose sides are $\overline{t}_1$, $\overline{x}_1$,$\cdots $,$\overline{x}_{n-2}$, $\overline{t}_{n-1}$.
This is a convex open set is   of ${\Bbb R}^{n-3}$.

On the other hand, the area function
$F_v\buildrel {\cal A} \over \rightarrow {\Bbb R}$ is given, for $t=(t_2,t_3,\cdots ,t_{n-2})\in F_v$, by:
$${\cal A}(t)
=h(\overline{t}_1,\overline{x}_1,t_2)+h(t_2,\overline{x}_2,t_3)+\cdots +h(t_{n-2},\overline{x}_{n-2},\overline{t}_{n-1}).\leqno{(5.3)}$$
Setting:
$$\cases{\displaystyle \omega
_1=(\overline{t}_1,\overline{x}_1,t_2)\cr {}\cr
\displaystyle
\omega _k=(t_k,\overline{x}_k,t_{k+1}) \hbox{ for $k\in \{2,\cdots ,n-3\}$} \cr {}\cr
\omega_{n-2}=(t_{n-2},\overline{x}_{n-2},\overline{t}_{n-1}),}
$$
one can express the partial derivatives of
${\cal A}$ as follows:
$${{\partial {\cal A}}\over {\partial t_k}}(t)={{\partial
}\over {\partial
t_k}}\left[h\left(t_{k-1},\overline{x}_{k-1},t_k\right)+h\left(t_k,\overline{x}_k,t_{k+1}\right)\right]=
{{\partial h}\over {\partial z}}(\omega _{k-1})+{{\partial
h}\over {\partial x}}(\omega _k).$$
A critical point $t\in F_v$ of the area function must satisfy:
$$
{{\partial h}\over {\partial z}}(\omega _{k-1})+{{\partial
h}\over {\partial x}}(\omega _k)=0 \; \; \hbox{and then} \; \;
\left({{\partial h}\over {\partial z}}(\omega_{k-1})\right)^2-\left({{\partial h}\over {\partial x}}(\omega_k)\right)^2=0.$$
Setting $\omega _{k-1}=(u,v,w)$  and  $\omega_k=(w,s,r)$, we obtain ${\alpha \over \beta }=0$ where:
$$\begin{array}{rcl}
{\alpha } = &w^2 \left(r^2uv + rsu^2 + rsv^2 - rsw^2 + s^2uv - uvw^2\right)\cr &\left(-r^2uv + rsu^2 + rsv^2 - rsw^2 - s^2uv + uvw^2\right).
\end{array}
$$
and $\beta = (r-s-w) (r+s-w) (r-s+w) (r+s+w)(-u-v+w) (u-v+w)
(-u+v+w) (u+v+w)$
or:
$$\begin{array}{rcl}
\left[(rsu^2 + rsv^2 - rsw^2) + (s^2uv - uvw^2 + r^2uv)\right]&\cr \left[(rsu^2 + rsv^2 - rsw^2) - (s^2uv - uvw^2 + r^2uv)\right]
= 0.
\end{array}
$$
This implies:
$$\left(rsu^2 + rsv^2 - rsw^2\right)^2 - \left(s^2uv - uvw^2 + r^2uv)\right)^2 = 0$$
or:
$$rs(u^2+v^2-w^2)=\pm uv(r^2+s^2-w^2).$$
Then we deduce:
$$\cos \widehat{M_{k-1}}=\pm \cos
\widehat{M_{k+1}}$$
which implies that the quadrilateral $(M_n,M_{k-1},M_k,M_{k+1})$ is inscribable for any
index $k\in \{2,\cdots , n-2\}$. Hence the  polygon $t$ is inscribable.

Conversely, if $t \in F_v$ is inscribable then, one can prove easily  that $t$ is a critical point of the area function.

\vskip0.8cm

\centerline{\tfp References}

\vskip 0.5cm

\noindent [Han] {\sc Hansen, V.L.} {\it Shadows of the Circle.} World Scientific (1998).

\smallskip

\noindent [Khi] {\sc Khimshiashvili, G.} {\it Cyclic polygons as critical points.} Proc. of I. Vekua Institute of Applied Mathematics Vol. 58, (2008), 74-83.

\smallskip

\noindent [Pen] {\sc Penner, R.C.} {\it The Decorated Teichm\"uller Space of Punctured Surfaces.} Commun. Math. Phys. 113, (1987), 299-339.

\smallskip

\noindent [Per] {\sc Perrin, D.} {\it  Math\'ematiques d'\'ecole. Nombres, mesures et g\'eom\'etrie.} Cassini, Paris, (2005).

\smallskip

\noindent [Ram] {\sc Ramon, P.} {\it Introduction \`a la g\'eom\'etrie diff\'erentielle discr\`ete.} Editions Ellipses, (2013).

\vskip1cm

\noindent {\pg Universit\'e Polytechnique Hauts-de-France

\noindent LAMAV, FR CNRS 2956,

\noindent  ISTV2, Le Mont Houy

\noindent 59313 Valenciennes Cedex 9

\noindent FRANCE}

\bigskip

\noindent{\pg aziz.elkacimi@uphf.fr}

\noindent{\pg abdellatif.zeggar@uphf.fr}

\end{document}